\input amstex
\documentstyle{amsppt}

\hsize=4.75in
\vsize=8in
\NoBlackBoxes
\def\today{\ifcase\month\or
 January\or February\or March\or April\or May\or June\or
 July\or August\or September\or October\or November\or December\fi
 \space\number\day, \number\year}

\def\nind{\noindent}
\def\RR{\text{{\rm I \hskip -5.75pt R}}}

\def\CC{\;\text{{\rm \vrule height 6pt width 1pt \hskip -4.5pt C}}}

\def\As{{\Cal A}}
\def\Bs{{\Cal B}}
\def\Cs{{\Cal C}}

\def\Hs{{\Cal H}}

\def\Ms{{\Cal M}}

\def\Os{{\Cal O}}

\def\Rs{{\Cal R}}
\def\Ss{{\Cal S}}
\def\Ts{{\Cal T}}

\def\net{{\{\As(\Os)\}_{\Os\in\Rs}}}

\def\idty{{\leavevmode\hbox{\rm 1\kern -.3em I}}}
\baselineskip 15pt plus 2pt
\spaceskip=.5em plus .25em minus .20em
\redefine\qed{\hbox{$\boxed{}$}}
\def\vN{von Neumann} 
\def\set#1{\left\lbrace#1\right\rbrace} 
\def\stt{\,\vrule\ }
\def\norm#1{\left\Vert#1\right\Vert} 
\def\abs#1{\left\vert#1\right\vert} 
\def\idty{{\leavevmode\hbox{\rm 1\kern -.3em I}}}
\def\RR{\text{{\rm I \hskip -5.75pt R}}}
 
\def\Bh{\Bs(\Hs)}  
\def\Aut{\text{\rm Aut}} 
\def\bracks#1{\lbrack#1\rbrack}
\redefine\qed{\hfill\break\rightline{\hbox{$\boxed{}$}}}
\def\Aoo{{\As(\Os_1)}}
\def\Aot{{\As(\Os_2)}}
\def\net{{\set{\As(\Os)}_{\Os\in\Rs}}}

\def\Awo{{\As(W_1)}}
\def\Awt{{\As(W_2)}}

\def\vacrept{{(\Hs,\pi,U(\RR^4),\Omega)}}

\def\pair{{(\As,\Bs)}}
\def\pairf{{(\As(\Os_1),\As(\Os_2))}}
\def\bellpas{{\beta(\phi,\As,\Bs)}}
\def\bell{{\beta(\hat{p},\As,\Bs)}}
\def\dual{(\hat{p},\As,\Bs)}
\def\viol{{\frac{1}{2}(A_1(B_1 + B_2) + A_2(B_1 - B_2))}}
\def\viol2{{(A_1(B_1 + B_2) + A_2(B_1 - B_2))}}
\def\Tvl{{\Cal T \pair }}              
\def\Tcl{\hbox{$\overline{\Ts}\pair$}} 
\def\ABs{\As,\Bs} 

\def\Poinc{{{\Cal P}^{\uparrow}_+}}
\rightheadtext {Bell's Inequalities and Algebraic Structure}
\leftheadtext {Stephen J. Summers}
\topmatter
\title
Bell's Inequalities and Algebraic Structure 
\footnotemark
\endtitle
\footnotetext{This paper is an extended version of a talk given at the 
International Conference on Operator Algebras and Quantum Field Theory, held
at the Accademia dei Lincei, Rome, in July, 1996.}
\author
Stephen J. Summers
\endauthor
\affil
Department of Mathematics \\ 
University of Florida \\
Gainesville, FL 32611, USA 
\endaffil
\abstract{We provide an overview of the connections between Bell's 
inequalities and algebraic structure.}
\endabstract
\date December 17, 1996 \enddate
\endtopmatter
\document

\heading 1. Introduction  \endheading

     Motivated by the desire to bring into the realm of testable hypotheses
at least some of the important matters concerning the interpretation of 
quantum mechanics evoked in the controversy surrounding the 
Einstein-Podolsky-Rosen paradox \cite{18}\cite{5}, Bell discovered the
first example \cite{3}\cite{4} of a family of inequalities
which are now generally called Bell's inequalities. These inequalities
provide an upper bound on the strength of correlations between systems
which are no longer interacting but have interacted in the past. Stated
briefly, Bell showed that if the correlation experiments can be modelled by
a single classical probability measure, then the strength of the correlations
must satisfy a bound which is violated by certain quantum mechanical
predictions (and, as has been verified experimentally, by nature - for a review
of this original application of Bell's inequalities and the experiments 
performed, see \cite{14}\cite{15}). Hence, it was established that in 
quantum mechanics and in nature there are correlations which cannot be modelled
by ``local hidden-variable theories''. Though there 
are many other applications of Bell's inequalities besides the one due to Bell 
(see {\it e.g.} \cite{26}\cite{42}\cite{40}), in this paper we shall 
concentrate on the relation between Bell's inequalities and algebraic 
structure which unexpectedly emerged in our study of Bell's inequalities in 
the context of quantum field theory.  \par
     First, we need to specify the mathematical context more precisely. As 
suggested by Ludwig's approach to statistical theories \cite{24} (see, in
particular, \cite{41}), the minimal amount of structure
required to model correlation experiments involving two subsystems (the only
situation we shall discuss here) is a so-called correlation duality
$(\hat{p},\As,\Bs)$, composed of two order unit spaces $\As$ and $\Bs$ (real
vector spaces with a vector ordering $\geq$ and a unit 1) and a bilinear
function $\hat{p} : \As \times \Bs \mapsto \RR$ such that $A \in \As$,
$B \in \Bs$, and $A,B \geq 0$ imply $\hat{p}(A,B) \geq 0$ and such that
$\hat{p}(1,1) = 1$. The function $\hat{p}$ models the preparation of the 
ensemble of systems. The observables of the
subsystem corresponding, for example, to $\As$ are represented by partitions 
$\{ A_i \mid i \in I \}$ of the unit in $\As$: $\sum_i A_i = 1$ with
$A_i \geq 0$ for each $i \in I$. Each $i \in I$ corresponds to a possible
outcome of the experiment. The probability (relative frequency) of the
joint occurrence of $i \in I$ and $j \in J$ in the two subsystems, 
respectively, is then given by $\hat{p}(A_i,B_j)$.  \par
     In quantum theory this structure is supplemented with additional
assumptions, of which we only mention in this introduction that the order unit
spaces $\As$ and $\Bs$ are posited to be the Hermitian part of subalgebras
(again denoted by $\As$ and $\Bs$) of a unital $C^*$-algebra $\Cs$ and that 
they commute elementwise. Moreover, the bilinear function $\hat{p}$ is given by
a state $\phi$ on $\Cs$: $\hat{p}(A,B) \equiv \phi(AB)$.  \par

     In \cite{32}\cite{34} we defined the {\it maximal Bell correlation} 
$\beta(\hat{p},\As,\Bs)$ in a correlation duality $(\hat{p},\As,\Bs)$ to be
$$\align
\beta(\hat{p},\As,\Bs) &\equiv \frac{1}{2}\sup\big(\hat{p}(A_1,B_1) +
\hat{p}(A_1,B_2) + \hat{p}(A_2,B_1) - \hat{p}(A_2,B_2)\big) \\
  &= \frac{1}{2}\sup\big(\hat{p}(A_1,B_1 + B_2) + \hat{p}(A_2,B_1 - B_2)\big)
\, , 
\endalign $$
where the supremum is taken over all elements $A_i \in \As$,
$B_i \in \Bs$ satisfying $-1 \leq A_i,B_i \leq 1$, $i = 1,2$. In the case
where $\As$ and $\Bs$ are the Hermitian parts of $C^*$-algebras 
$\As, \Bs \subset \Cs$ and $\hat{p}$ is given by a state $\phi$, then we shall
write instead
$$\beta(\phi,\As,\Bs) = \frac{1}{2}\sup\phi(A_1(B_1 + B_2) + A_2(B_1 - B_2))
\, . $$

     In \cite{34}, the following bounds for $\bell$ were proven.

\proclaim{Theorem 1.1} \cite{32}\cite{34} (a) If $\dual$ is an arbitrary 
correlation duality, then 
$$\bell \leq 2 \, . \tag{1.1} $$

   (b) If $\As$ (or $\Bs$) is the Hermitian part of some $C^*$-algebra,
then 
$$\bell \leq \sqrt{2} \, , \tag{1.2} $$
for every admissible $\hat{p}$ as described above. \par
   (c) If $\As$ (or $\Bs$) is the Hermitian part of some abelian $C^*$-algebra,
then
$$\bell \leq 1 \, , \tag{1.3} $$
for every admissible $\hat{p}$ as described above. 
\endproclaim

     The estimate (1.2) was first given in \cite{12} and subsequently 
rediscovered by a number of researchers. In fact, it has been shown 
\cite{23} that if $\As$ and $\Bs$ are merely distributive real algebras with
identity, 
then the inequality $\bellpas \leq \sqrt{2}$ must hold for any state $\phi$. 
Hence, the bound (1.2) is also satisfied by Jordan algebras and distributive 
Segal algebras. It is easy to construct correlation dualities $\dual$ 
saturating the bound (1.1), but it can also be saturated for suitable choice 
of $\dual$ with $\As$ and $\Bs$ nondistributive real algebras (see {\it e.g.} 
\cite{22}). We shall see below that the bound (1.2) can also be
saturated. \par
     Initially, the purpose of these inequalities was to serve as necessary
conditions for their respective hypotheses. In particular, considering (1.2), 
a correlation experiment in the laboratory yielding 
$$\hat{p}(A_1,B_1) + \hat{p}(A_1,B_2) + \hat{p}(A_2,B_1) - \hat{p}(A_2,B_2) 
> \sqrt{2} $$
could not thus be described by quantum theory, as pointed out by Cirel'son
\cite{12}.\footnote{This has not been observed in the laboratory.} 
The bound (1.3) is equivalent \cite{34} to the CHSH version \cite{13} of 
Bell's inequality. There are many different proofs of 
this inequality in the literature - at least one for each metatheoretical 
framework for the discussion of such correlation experiments. In \cite{33} 
we proved a stronger result than that given here, namely, that if the 
correlation duality is quasi-classical\footnote{See \cite{41} for a 
definition of quasi-classical correlation duality.}, then (1.3) must hold for 
every admissible $\hat{p}$. But for our present purposes, the special case 
indicated above will suffice, since it makes clear that if the given 
correlation experiment can be modelled by a classical theory (thereby yielding
abelian observable algebras), then Bell's inequality (1.3) must hold. Any 
correlation experiment yielding
$$\hat{p}(A_1,B_1) + \hat{p}(A_1,B_2) + \hat{p}(A_2,B_1) - \hat{p}(A_2,B_2) 
> 1 $$
could therefore not be modelled by a classical theory (at least a 
classical theory providing a correlation duality). And since quantum theory 
predicts (see below) and nature confirms the existence of such correlations
violating Bell's inequality, the reader may begin to grasp the significance of
Bell's original intention. \par
     Henceforth, we shall restrict our attention to the case where $\As$ and 
$\Bs$ are (Hermitian parts of) $C^*$-subalgebras of a unital $C^*$-algebra
$\Cs$, each containing the identity $\idty \in \Cs$. In Section 2 we shall 
examine the properties of the maximal Bell correlation $\beta(\phi,\As,\Bs)$
and shall use it to define an algebraic invariant $\beta(\As,\Bs)$ of the
(isomorphism class of the) pair $\pair$. If $\As$ and $\Bs$ are realized as a
pair of commuting von Neumann algebras on a separable Hilbert space $\Hs$, 
it would follow that $\As$ is contained in $\Bs'$, the commutant of $\Bs$. 
Therefore, defining $\beta\pair$ for $\pair$ is equivalent to defining an 
algebraic invariant of the (isomorphism class of the) inclusion 
$\As \subset \Bs'$ of \vN\ algebras - an invariant which is quite distinct 
from that of Jones and others (see, {\it e.g.} \cite{27}). We 
shall also explain in Section 2 which algebraic structural properties are 
associated with the maximal possible value of the invariant $\beta\pair$. 
There were some surprises here. In Section 3 we shall report on results 
obtained in the special context of algebraic quantum field theory, which, 
among other things, establish that $\beta\pair$ takes on infinitely many 
distinct values for suitable choices of $\pair$.  \par
 
\bigpagebreak 

\heading 2. Maximal Bell Correlations and Algebraic Invariants \endheading

     As shall be made clear later, it is most useful to assume that the algebra
$\Cs$ above is the algebra $\Bs(\Hs)$ of all linear bounded operators on a 
separable Hilbert space $\Hs$, and that $\As$ and $\Bs$ are commuting \vN\
subalgebras of $\Bs(\Hs)$. In order to state assertions succinctly, we set for 
such algebras $\ABs\subset\Bh$  
$$\align
 \Tvl \equiv \Bigl\lbrace
           {1\over2}\bigl(A_1(B_1+B_2)&+A_2(B_1-B_2)\bigr)
                  \Big\vert
              \\&A_i = A_i^* \in\As,\ B_i = B_i^* \in\Bs,\  
                       -\idty\leq A_i,B_i \leq\idty 
                      \Bigr\rbrace \, .
\endalign$$
The set $\Tvl$ contains all the observables relevant for 
testing violations of Bell's inequalities in the (independent) 
subsystems modelled by the (mutually commuting) pair of algebras of 
observables $\pair$. The $\sigma\bigl(\Bh,\Bh_*\bigr)$-closed convex 
hull of $\Tvl$ will be denoted by $\Tcl$. Since in such 
expressions we can replace $A_i$ by $\lambda A_i$ with 
$-1\leq\lambda\leq1$, the image 
$\phi\bigl(\Tvl\bigr)\subset\RR$ of $\Tvl$ under any 
state $\phi\in\Bh^*$ is a symmetric interval 
\footnote{In particular, $\phi(\Tvl)$ is a connected set.} 
around zero. This interval is therefore characterized by the 
maximal Bell correlation  
$$ \beta(\phi,\ABs) = \sup\set{ \phi(T)\stt T\in\Tvl} 
\, . $$
Of course, if the state $\phi$ is normal on $\Bh$, then we also have
$$ \beta(\phi,\ABs) = \sup\set{ \phi(T)\stt T\in\Tcl} 
\, .$$
Since $\idty\in\Tvl$ (set $A_i=B_i=\idty$), we have 
$\beta(\phi,\ABs)\geq1$. As mentioned above, any value larger than $1$ is a 
{\it violation} of Bell's inequalities in the given state, detectable by 
observables $A_i,B_i$ in the subalgebras $\ABs$. On the other hand, by Theorem
1.1, the largest possible value for $\beta(\phi,\ABs)$ in any state is 
$\sqrt2$, and if this value is attained we speak of a {\it maximal violation} 
of Bell's inequalities by the pair of algebras $\pair$ in the state $\phi$. 

     We summarize some fundamental properties of maximal Bell correlations in 
the following lemma. 

\proclaim{Lemma 2.1} \cite{34}\cite{39} Let $\ABs\subset\Bh$ be mutually 
commuting \vN\ algebras and $\phi \in \Bh^*$ be a state. Then all of the 
following assertions are true. 
\roster
\item"1." $1\leq\beta(\phi,\ABs)\leq\sqrt2$. 

\item"2." The map $\phi\mapsto\beta(\phi,\ABs)$ is convex. 

\item"3." 
The map $\phi\mapsto\beta(\phi,\ABs)$ is lower semicontinuous in the 
$\sigma(\Bh^*,\Bh)$ topology. 

\item"4." 
$\abs{\beta(\phi,\ABs)-\beta(\psi,\ABs)}\leq\sqrt{2}\norm{\phi-\psi}$, 
so, in particular, $\beta(\phi,\ABs)$ is norm continuous in the 
state $\phi$. 

\item"5." 
If either $\As$ or $\Bs$ is abelian, then $\beta(\phi,\ABs)=1$ for 
all states $\phi \in \Bh^*$. If both of the algebras are nonabelian, 
there exists a normal state $\phi \in \Bh_*$ such that 
$\beta(\phi,\ABs)=\sqrt2$. 

\item"6." 
If $\phi$ is a convex combination of product states, then 
$\beta(\phi,\ABs)=1$. 

\item"7." 
There exists a (possibly non-normal) state $\phi$ with $\beta(\phi,\ABs)=1$. 
\endroster 
\endproclaim
     
     Note that Lemma 2.1 (5) asserts that for {\it any} choice of a pair
$\pair$ of nonabelian, mutually commuting \vN\ algebras, there exists a normal
state $\phi$ such that the upper bound in (1.2) is attained, {\it i.e.} 
Bell's inequality is maximally violated. As explained in \cite{38},
though the basic idea goes all the way back to Bell, if two commuting algebras
$\As$ and $\Bs$ contain copies of the set $M_2(\CC)$ of two-by-two complex 
matrices, they maximally violate Bell's inequality in some normal state.
And such copies of $M_2(\CC)$ can always be found in nonabelian algebras.
It is interesting that the presence of copies of $M_2(\CC)$ is also
{\it necessary} for maximal violation of Bell's inequalities. 

\proclaim{Proposition 2.2} \cite{34} Let $\pair$ be a pair of commuting
subalgebras of a unital $C^*$-algebra $\Cs$. If $A_i \in \As$, $B_j \in \Bs$,
are selfadjoint and bounded in norm by $1$, and if $\phi$ is a state
on $\Cs$ with its restrictions to both $\As$ and $\Bs$ faithful, then the
equality 
$$\frac{1}{2}\phi(\viol2) = \sqrt{2}$$
entails that $A_i^2 = \idty$, $i = 1,2$, and $A_1A_2 + A_2A_1 = 0$
(similarly for the $B_j$). Therefore, $A_1,A_2$ and 
$A_3 \equiv -\frac{i}{2}[A_1,A_2]$ form a realization of the Pauli spin 
matrices in $\As$ and hence generate a copy of $M_2(\CC)$ in $\As$ (similarly 
for the $B_j$ in $\Bs$). Moreover, the $A_i$, resp. the $B_j$, are contained 
in the centralizer of $\As$ in $\phi$, resp. centralizer of $\Bs$ in $\phi$. 
\endproclaim

\nind This has consequences for experimental physics: if one wishes to 
measure a maximal violation of Bell's inequalities, one must observe
quantities which can be modelled by Pauli spin matrices. \par 

     By Lemma 2.1 (4) and the norm-continuity of the map 
$\lambda\mapsto\lambda\phi+(1-\lambda)\psi$, it is clear that the 
range of the convex functional $\beta(\cdot,\ABs)$ on the state 
space is an interval. According to Lemma 2.1 (5), this interval 
contains $\sqrt2$ in all cases of interest, {\it i.e.} where neither 
algebra is abelian. Considering the range of $\beta(\cdot,\ABs)$ 
over {\it all} states $\phi \in \Bs(\Hs)^*$, we see from Lemma 2.1 
(7) that the lower endpoint will always be attained and be equal to 
$1$. However, restricted to the {\it normal} state space the infimum 
may be strictly between $1$ and $\sqrt2$ (and this indeed occurs - 
see below). We summarize such facts in the next lemma. 

\proclaim{Lemma 2.3} \cite{39} Let $\ABs\subset\Bh$ be commuting 
nonabelian von Neumann algebras. 
\roster
\item"1." 
The range of the map $\phi\rightarrow\beta(\phi,\ABs)$ as $\phi$ runs
through the set of all states on $\Bh$ is the closed interval 
$\bracks{1,\sqrt{2}}$. 

\item"2."
The range of the map $\phi\rightarrow\beta(\phi,\ABs)$ as $\phi$ 
runs through the set of all normal states on $\Bh$ is the interval 
$\lbrace c,\sqrt{2}\rbrack\,$, 
\footnote{The interval may or may not contain the point $c$, and it 
may be degenerate, {\it i.e.} equal to the singleton set $\set{\sqrt{2}}$.} 
where $c \in \bracks{1,\sqrt{2}}$. 

\item"3." 
If $c < \sqrt{2}$, then there exists a norm dense set of normal states
$\psi$ such that $\beta(\psi,\As,\Bs) < \sqrt{2}$. 
\endroster
\endproclaim

     The number $c$ appearing in Lemma 2.3 is an invariant of the (isomorphism 
class of the unordered) pair $\pair$ of \vN\ algebras. We therefore 
define, for any pair $(\ABs)$ of commuting \vN\ algebras on a 
Hilbert space $\Hs$, the {\it Bell correlation invariant} of the 
pair to be this number: 
$$\beta(\ABs) \equiv 
    \inf\set{\beta(\phi,\ABs) \stt 
          \phi \,\text{ a normal state on} \quad \Bh } 
\, .$$
A variant of this which will be of interest is 
$$ \beta_*(\ABs)= \sup\set{\lambda\in\RR\stt 
                  \exists T\in\Tcl\ \text{with }\ 
                          T\geq\lambda\idty}
\, .$$
Thus $\beta(\ABs)\geq\lambda$ means that for any normal state 
$\phi$ there is an operator $T\in\Tcl$ (or, equivalently, an 
operator $T\in\Tvl$) with $\phi(T)\geq\lambda$. If this operator can 
be chosen independently of $\phi$, then we have 
$\beta_*(\ABs)\geq\lambda$. From this discussion it is clear that 
$$ \beta_*(\ABs)\leq\beta(\ABs)\leq\beta(\phi,\ABs) 
\, , $$
for all normal states $\phi \in \Bh_*$. We emphasize that a 
consequence of Lemma 2.3 is that if $\pair$ is any pair of mutually 
commuting, nonabelian von Neumann algebras, then for {\it any} value 
$r \in (\beta\pair,\sqrt{2}\rbrack$ there exists a normal state $\phi \in 
(\As \vee \Bs)_*$ such that $\beta(\phi,\As,\Bs) = r$, where $\As \vee \Bs$ 
denotes the subalgebra of $\Bh$ generated by $\As$ and $\Bs$. The following 
result was communicated to us by Shulman, and it establishes that 
$\beta\pair$ and $\beta_* \pair$ are alternative ways to compute the same 
algebraic invariant.

\proclaim{Proposition 2.4} \cite{30} For any pair $\pair$ of commuting \vN\ 
algebras, one has $\beta(\As,\Bs) = \beta_*(\As,\Bs)$.
\endproclaim

\demo{Proof} It is easy to see that the definition entails
$$\beta_*(\As,\Bs) = \sup\{ \lambda(T) \mid T \in \Tcl \} \, , $$
where 
$\lambda(T) \equiv \sup\{\lambda\in\RR \mid T \geq \lambda\cdot 1 \}$. It 
follows that 
$$\lambda(T) = \inf\{\phi(T) \mid \phi \in \Ss(\Bs(\Hs)) \} = 
\inf\{\phi(T) \mid \phi \in \Ss_*(\Bs(\Hs)) \} \, , $$
where $\Ss(\Bs(\Hs))$, resp. $\Ss_*(\Bs(\Hs))$, is the set of all 
states, resp. normal states, on $\Bs(\Hs)$.
Hence, since $\Tcl$ is *-weakly compact and the function 
$(\phi,T) \mapsto \phi(T)$ is continuous on $\Ss(\Bs(\Hs)) \times \Tcl$, one
has 
$$\beta_*\pair = \sup_{T\in\Tcl}\inf_{\phi\in\Ss_*(\Bs(\Hs))}\phi(T)
=  \inf_{\phi\in\Ss_*(\Bs(\Hs))}\sup_{T\in\Tcl}\phi(T) \, ,$$
using the generalization of Ky Fan's minimax result in Prop. 1 of \cite{6}. 
Since $\sup_{T \in \Tcl}\phi(T) = \sup_{T \in \Tvl}\phi(T)$, one finds
$$\beta_*\pair = \inf_{\phi\in\Ss_*(\Bs(\Hs))}\bellpas = \beta\pair \, , $$
as asserted.
\hfill\qed\enddemo

\nind Hence we may conclude from the comment above that if 
$\beta\pair \geq \lambda$, then there exists an element $T \in \Tcl$ such that 
$\phi(T) \geq \lambda$ for {\it all} normal states $\phi$. (For the maximal
possible value, $\lambda = \sqrt{2}$, this was established earlier 
\cite{36}.) \par

     Some basic properties of this invariant are given next. 

\proclaim{Proposition 2.5} \cite{39} Let $\ABs\subset\Bh$ be commuting 
\vN\ algebras. 
\roster
\item"1." 
If $\As_1\subset\As_2$ and $\Bs_1\subset\Bs_2$, then 
$\beta(\As_1,\Bs_1)\leq\beta(\As_2,\Bs_2) \, .$            

\item"2." 
Let $\As=\underset{i \in I}\to{\oplus}\As_i$, and 
$\Bs=\underset{i \in I}\to{\oplus}\Bs_i$, with $I$ an arbitrary 
index set. Then 
$\beta(\As,\Bs)=\hbox{\rm inf}_{i \in I}\beta(\As_i,\Bs_i) \, .$ 

\item"3." 
If the pair $\pair$ is split, {\it i.e.} if there is a type $I$ factor
$\Ms$ with $\As\subset\Ms\subset\Bs'$, then $\beta(\ABs)=1$ . 

\item"4."
If $\beta\pair > 1$, then both $\As$ and $\Bs$ are nonabelian. Moreover,
if $\beta\pair > 1$, there exist no normal states in the weak*-closure of
the convex hull of the product states across $\pair$, {\it i.e.} all
normal states on $\As\vee\Bs$ are ``entangled''.
 
\endroster
\endproclaim

    In the terminology of \cite{37}, if $\beta\pair$ takes the value
$\sqrt{2}$, we say that the pair $\pair$ is {\it maximally correlated}, since 
$\beta\pair = \sqrt{2}$ entails that Bell's inequalities are {\it maximally} 
violated in {\it every} normal state across $\pair$. A very useful fact is that
a pair $\pair$ is maximally correlated if and only if there is maximal 
violation of Bell's inequalities in at least one faithful normal state on 
$\As \vee \Bs$.

\proclaim{Proposition 2.6} \cite{37} Let $\pair$ be a pair of commuting
von Neumann algebras acting on a separable Hilbert space $\Hs$. Then the
following are equivalent.
\roster
\item"1." $\pair$ is maximally correlated.

\item"2." There exists a faithful state $\phi \in (\As \vee \Bs)_*$ such
that $\bellpas = \sqrt2$.
\endroster
\endproclaim

\nind Thus, in order to check if a pair is maximally correlated, it suffices
to check the value of $\bellpas$ in any conveniently chosen faithful normal 
state $\phi$. A structural characterization of maximally correlated pairs of 
\vN\ algebras was given in \cite{37}.

\proclaim{Theorem 2.7} \cite{37} Let $\pair$ be a pair of commuting \vN\
algebras acting on a separable Hilbert space $\Hs$. Then the following are 
equivalent. 
\roster
\item"1." The pair $\pair$ is maximally correlated.

\item"2." There exists a type $I$ factor $\Ms \subset \As \vee \Bs$ such
that $\As \cap \Ms$ and $\Bs \cap \Ms$ are spatially isomorphic to the 
hyperfinite type $II_1$-factor $\Rs_1$ and are relative commutants of each 
other in $\Ms$.

\item"3." Up to unitary transformation, $\Hs = \Hs_1 \otimes \tilde{\Hs}$, 
where $\Hs_1$ is the standard representation space of $\Rs_1$, and there exist 
von Neumann algebras 
$\tilde{\As} \subset \tilde{\Bs}' \subset \Bs(\tilde{\Hs})$ such that 
$\As = \Rs_1 \otimes \tilde{\As}$ and $\Bs = \Rs_1 ' \otimes \tilde{\Bs}$.

\endroster
\endproclaim

     In the special case $\Bs = \As'$, we set 
$\beta(\As) \equiv \beta(\As,\As')$, which yields a new invariant of the \vN\ 
algebra $\As$ itself. Surprisingly, in this special case the property that
$(\As,\As')$ is maximally correlated ({\it i.e.} $\beta\pair = \sqrt2$) can be 
restated in terms of the algebraic properties $L_{\lambda}$, resp. 
$L_{\lambda}'$, which were isolated and studied some 25 years ago by Powers 
\cite{28}\cite{29}, resp. by Araki \cite{1}, as part of the program to 
classify von Neumann algebras. For the reader's convenience, we first recall 
the definitions of these properties. \par
     Let $\As$ be a $C^*$-algebra with unit $\idty$. Then $N \in \As$ is 
called a $I_2$-generator if $N^2 = 0$ and $NN^* + N^*N =  \idty$. Let 
$V_{\As}$ denote the set of $I_2$-generators in $\As$. Clearly, if $N$ is 
contained in $V_{\As}$, then $N^*N$ and $NN^*$ are nonzero complementary 
projections, {\it i.e.} their sum is $\idty$ and their product is zero. 
Moreover, the $C^*$-algebra generated by $N$ is isomorphic to $M_2(\CC)$ and 
contains the unit $\idty$ of $\As$. Conversely, if $\As$ contains a copy of 
$M_2(\CC)$ containing $\idty$, then $V_{\As} \neq \emptyset$. Note that if 
$A_i \in \As$ satisfies $A_i^* = A_i$, $A_i^2 = \idty$ and 
$A_1A_2 + A_2A_1 = 0$ (which, by Prop. 2.2, is the case if $A_1,A_2$ are 
maximal violators of Bell's inequalities in some faithful state on $\As$), 
then $N \equiv \frac{1}{2}(A_1 + iA_2)$ is an element of $V_{\As}$. 

\proclaim{Definition}\cite{29}\cite{1} A von Neumann 
algebra $\As$ is said to have 
property $L_{\lambda}$ (resp. $L_{\lambda}'$) with $\lambda \in [0,1/2]$ if
for every $\epsilon > 0$ and any normal state $\phi$ on $\As$ (resp.
finite family $\{\phi_i\}_{i=1}^n$ of normal states on $\As$), there exists
an $N \in V_{\As}$ such that for any $A \in \As$,
$$\vert\lambda\phi(AN) - (1-\lambda)\phi(NA)\vert \leq \epsilon\Vert A\Vert
           \tag{2.1}$$
(resp. for any $A \in \As$ and $i=1,\ldots,n$
$$\vert\lambda\phi_i(AN)-(1-\lambda)\phi_i(NA)\vert \leq 
                                          \epsilon\Vert A\Vert) \, . $$
\endproclaim

     An alternative characterization of property $L_{\lambda}'$ has been made 
in terms of the asymptotic ratio set \cite{2}\cite{1} of the algebra $\As$. 
The asymptotic ratio set $r_{\infty}(\As)$ of a von Neumann algebra $\As$ is 
the set of all $\alpha \in \bracks{0,1}$ such that $\As$ is $W^*$-isomorphic 
to $\As \overline{\otimes} \Rs_{\alpha}$, where 
$\set{\Rs_{\alpha}}_{\alpha \in \bracks{0,1}}$ is the family of hyperfinite 
factors constructed by von Neumann \cite{25}, Powers \cite{28} and Araki and 
Woods \cite{2}. It is known that property $L_{\lambda}'$ is strictly stronger 
than property $L_{\lambda}$ \cite{1}, that property $L_{\lambda}'$ implies
property $L_{1/2}'$ for any $\lambda$ \cite{1}\cite{2}, and that 
property $L_{\lambda}'$
for $\As$ is equivalent to $\frac{\lambda}{1-\lambda} \in r_{\infty}(\As)$ 
\cite{1}. Using Prop. 2.2 one easily sees that if $A_1,A_2 \in \As$ 
are maximal violators of Bell's inequalities in the normal state $\phi$ on
$\As \vee \Bs$, where $\Bs \subset \As '$, then $N \equiv 
\frac{1}{2}(A_1 + iA_2) \in V_{\As}$ satisfies (2.1) with $\epsilon = 0$
and $\lambda = 1/2$. We can now state our structure result for maximally
correlated pairs $(\As,\As')$. \par

\proclaim{Theorem 2.8} Let $\As$ be a \vN\ algebra in a separable 
Hilbert space $\Hs$ with cyclic and separating vector. Then the following 
conditions are equivalent \cite{36}. 
\roster
\item"1." 
$\beta(\As) = \sqrt{2} \, .$ 

\item"2." 
There exists a sequence $\set{T_n} \subset \Ts(\As,\As')$ such that
$T_n \rightarrow \sqrt{2}\cdot\idty$ in the $\sigma$-weak topology, {\it i.e.}
$\beta_*(\As) = \sqrt{2} \, .$ 

\item"3." 
$\As \simeq \As\otimes\Rs_1$, where $\Rs_1$ is the hyperfinite type 
$II_1$ factor, {\it i.e.} $\As$ is strongly stable. 

\item"4." 
$\As$ has property $L_{1/2}' \, .$ 
\endroster
Moreover, the following are equivalent to each other \cite{37}. 
\roster
\item"5." 
$\As$ has property $L_{1/2} \, .$ 

\item"6." 
For any vector state $\omega(A) = \langle\Omega,A\Omega\rangle$,
$\Omega \in \Hs$, one has $\beta(\omega,\As,\As') = \sqrt{2} \, .$
\endroster
\endproclaim

\nind Note that this result and Prop. 2.5 (1) entail that if $\pair$ is
maximally correlated, then both $\As$ and $\Bs$ are strongly stable and
have property $L_{1/2}'$. \par 
     On the one hand, there are many algebras having property 
$L_{1/2}'$, and, on the other, if $\As$ is a type $I$ factor, then 
$\beta(\As) = 1$, by Proposition 2.5 (3). Though one therefore has examples 
of von Neumann algebras with $\beta(\As)$ taking values at the endpoints
of the admissible interval $[1,\sqrt2]$, at present it is not known 
whether there exist factors such that $1 < \beta(\As) < \sqrt{2}$. 
It would be interesting to find such factors, since Theorem 2.8 would entail
that $\As$ could not have property $L_{\lambda}'$ for any $\lambda$, {\it i.e.}
$\beta(\As)$ would then be established as a nontrivial invariant of \vN\ 
algebras extending the family $L_{\lambda}'$. However, in the more 
general setting where $\Bs \neq \As'$, we shall see that there are 
indeed examples of pairs $\pair$ such that $1 < \beta\pair < \sqrt{2}$.

\bigpagebreak

\heading 3. Quantum Field Theory and the Bell Correlation Invariant 
\endheading

     The original motivation of this work was to study Bell's inequalities 
in the context of quantum field theory. But the results in this special case 
(as viewed from the standpoint of operator algebra theory) can also be used to 
show that the invariant $\beta\pair$ defined above can attain values other 
than $1$ and $\sqrt{2}$. We shall therefore specialize now to the context of 
algebraic quantum field theory (see {\it e.g.} \cite{21}).

     The basic object in AQFT is a net of $C^*$-algebras $\net$ associating to 
open subsets\footnote{The index set $\Rs$ can be a proper directed subset of 
the collection of all open subsets of Minkowski space.} $\Os$ of Minkowski 
space unital $C^*$-algebras $\As(\Os)$ satisfying the following 
basic axioms, which are naturally motivated by the interpretation of 
$\As(\Os)$ as the algebra generated by all observables which can be measured 
in the spacetime region $\Os$:

   1. {\it Isotony\/}: If $\Os_1 \subset \Os_2$, then 
$\As(\Os_1) \subset \As(\Os_2)$; hence, the inductive limit 
$C^*$-algebra $\As \equiv \underset{\Os \in \Rs}\to\vee 
\As(\Os)$ exists.

   2. {\it Locality\/}: If $\Os_1$ is spacelike separated from 
$\Os_2$, then every element of $\As(\Os_1)$ commutes with every 
element of $\As(\Os_2)$.

   3. {\it Poincar\'e Covariance\/}: There is a faithful 
representation 
$\Poinc \ni \lambda \mapsto \alpha_{\lambda} 
     \in \Aut\As$ of the identity component $\Poinc$ 
of the Poincar\'e group in the automorphism group $\Aut\As$ of 
$\As$ such that $\alpha_{\lambda}(\As(\Os)) = \As(\Os_{\lambda})$, 
where $\Os_{\lambda}$ is the image under $\lambda$ of the region 
$\Os$. 

     In addition, the set of interesting states (or representations) on the net
is selected by some further general considerations, for example the requirement
that the representation be Poincar\'e, or at least translation, covariant. 
Among these are vacuum representations $\vacrept$. 
In particular, in the representation space $\Hs$ there exist a unit vector 
$\Omega$ with corresponding vector state $\omega_0$ on $\As$ and a strongly 
continuous unitary representation $U(\RR^4)$ of the translation 
subgroup of the Poincar\'e group, whose joint generators satisfy the 
spectrum condition, which leaves $\Omega$ invariant, and which 
satisfies\footnote{For a purely algebraic characterization of vacuum states on 
Minkowski space, see \cite{11}.}
$$U(x)\pi(A)U(x)^{-1} = \pi(\alpha_x(A)), \qquad 
      \text{for all} \, x \in \RR^4 \, , \, A \in \As \, .
$$
When restricting one's attention to a single representation, it is 
convenient and customary to identify the algebras $\As(\Os)$ with the von 
Neumann algebras $\pi(\As(\Os))''$ on the representation space $\Hs$. \par
     A striking fact which emerged in our study of Bell's inequalities is that
in quantum field theory, as opposed to nonrelativistic quantum theory -
where the observable algebras for subsystems are typically type $I$, so
Prop. 2.5 (3) yields $\beta\pair = 1$ - there are many spacetime regions with
observable algebras evincing maximal violation of Bell's inequalities in every 
normal state, {\it i.e.} maximal violation independent of how the system has
been prepared. In the following, wedge regions in Minkowski space are
Poincar\'e transforms of the unbounded region 
$\{x \in \RR^4 \mid x^1 > \vert x^0 \vert \}$, and double cones are the 
nonempty interior of the bounded intersection of the forward, resp. 
backward, lightcone of two timelike separated points. For a given open set
$\Os \subset \RR^4$, $\Os'$ denotes the interior of the set of all points in 
$\RR^4$ which are spacelike separated from all points of $\Os$. \par

\proclaim{Theorem 3.1} 
\cite{32}\cite{34}\cite{35}\cite{36}\cite{37} For any net
$\net$ of local algebras satisfying the above assumptions the following
is true.
\roster
\item"1." In any vacuum representation, in any superselection sector
occurring in the Doplicher, Haag, Roberts theory of superselection structure
\cite{16}, and in any massive particle representation as described by 
Buchholz and Fredenhagen \cite{8}, one has $\beta(\As(W)) = \sqrt{2}$,
for all wedge regions $W \subset \RR^4$. Hence, if wedge duality is satisfied
({\it i.e.} if $\As(W)' = \As(W')$ for all wedges $W$), then $(\As(W),\As(W')$
is maximally correlated.

\item"2." In any dilatation-invariant vacuum representation satisfying wedge
duality, the pair $(\As(W_1),\As(W_2))$ is maximally correlated for any
spacelike separated wedges $W_1,W_2$, independent of the distance of 
separation.

\item"3." In any free field theory, in any locally Fock field theory 
({\it e.g.} $P(\phi)_2$, $Yukawa_2$, {\it etc.}), and in any 
dilatation-invariant theory satisfying wedge duality, the pair $\pairf$ is 
maximally correlated for any pair $(\Os_1,\Os_2)$ of tangent double cones in 
Minkowski space, hence $\beta(\As(\Os)) = \sqrt2$ for any double cone $\Os$.
\endroster
\endproclaim

\nind One sees that quantum field theory predicts maximal violation of Bell's
inequalities, quite independently of dynamics or even
preparation. We remark that, although such tangent algebras are maximally
correlated and therefore are not split, they still manifest very strong
independence properties - in particular, arbitrary normal states on the
subalgebras $\As(\Os_1)$, $\As(\Os_2)$ have simultaneous normal extensions to
$\As(\Os_1) \vee \As(\Os_2)$ \cite{19}.  \par
     Going beyond the cases mentioned in Theorem 3.1, it is known that in 
many concrete quantum field models \cite{7}\cite{17}\cite{31} (or under 
general, physically motivated assumptions \cite{9}\cite{10}), pairs of 
algebras $\pairf$ associated with strictly spacelike separated 
double cones are split. From Lemma 2.3 (2) we may conclude that under such 
circumstances, for any $r \in \bracks{1,\sqrt{2}}$ there exists a normal state 
$\phi$ such that $\beta(\phi,\As(\Os_1),\As(\Os_2)) = r$. Thus, in principle, 
every possible degree of violation of Bell's inequalities can be attained by 
suitable preparation on strictly spacelike separated double cone algebras, no 
matter how far apart their localizations are separated. However, in such 
cases one also has $\beta(\As(\Os_1),\As(\Os_2)) = 1$ (Prop. 2.5 (3)). On the 
other hand, pairs of algebras $\pairf$ associated with spacelike separated 
wedges are not split \cite{7}, allowing the {\it a priori} possibility that 
$1 < \beta(\As(\Os_1),\As(\Os_2)) < \sqrt{2}$ for such algebras, a possibility
which we shall see realized below.  

     In \cite{34} it was shown that in any irreducible vacuum sector 
with a strictly positive mass gap $m>0$ the maximal Bell correlation 
in the vacuum state $\omega_0$ satisfies an upper bound which 
decreases exponentially down to $1$:
$$\beta(\omega_0,\Aoo,\Aot) \leq 1 + 2e^{-md(\Os_1,\Os_2)} \, , 
\tag{3.1}$$
where $d(\Os_1,\Os_2)$ is the maximal timelike distance $\Os_{1}$ 
can be translated before $\Os_{1} \not\subset \Os_2 '$. Note that, in general, 
this is not the same as the smallest spacelike distance $d'(\Os_1,\Os_2) 
\equiv \inf\sqrt{-(x_1-x_2)^2}$ between the regions, although for 
wedges these two quantities coincide. The estimate (3.1) is a 
useful bound for large $d(\Os_1,\Os_2)$, but is clearly too crude 
for small distances 
\footnote{Note that it would be very interesting to obtain a lower 
bound on the quantity 
\newline $\beta(\omega_0,\Aoo,\Aot)$ which 
would decrease to 1 exponentially (with the same exponent $m$). This 
would entail that $\beta(\omega_0,\Aoo,\Aot)$ contains information 
about the lowest mass in the theory as well as metric information 
about the underlying space-time.}. 
The refinement of the bound (3.1) for small distances is based on the 
following proposition.  

\proclaim{Proposition 3.2} \cite{39} Let $\omega$ be a state, defined on 
an algebra containing two commuting C*-algebras, $\As$ and $\Bs$, and 
satisfying 
$$ \left\vert \omega(AB)-\omega(A)\omega(B) \right\vert
   \leq \Gamma\cdot \bigl(\omega(A^*A)\omega(AA^*)
                     \omega(B^*B)\omega(BB^*)\bigr)^{1/4} \, ,
$$
for some constant $0\leq\Gamma\leq1$, and all $A\in\As$, $B\in\Bs$, 
with $\Vert A \Vert \leq 1$, $\Vert B \Vert \leq 1$. Then
$$ \beta(\omega,\As,\Bs)
     \leq \min\Big\lbrace
             \sqrt2\  \frac{6+4\sqrt2+\Gamma}{7+4\sqrt2}\ ,\ 
             1+\sqrt{2}\,\Gamma \big\rbrace \, .
$$
\endproclaim

     Combined with the cluster bound available in vacuum representations with 
a mass gap \cite{20}, Prop. 3.2 yields the following result, which 
provides sharp short-distance bounds on the vacuum Bell correlation. 

\proclaim{Corollary 3.3} \cite{39} Let $\net$ be a local net in an 
irreducible vacuum representation $\vacrept$ with a mass gap 
$m>0$. Then, for any pair $(\Os_1,\Os_2)$  of spacelike separated 
regions, 
$$\beta(\omega_0,\As(\Os_1),\As(\Os_2)) \leq 
   \sqrt{2} - \frac{\sqrt{2}}{7+4\sqrt{2}}(1-e^{-md(\Os_1,\Os_2)}) 
\, .$$
\endproclaim

    Therefore, in a vacuum representation with a mass gap and
for any two convex spacetime regions $\Os_1$ and $\Os_2$ with a nonzero 
spacelike distance between them, Corollary 3.3 yields the bound 
$\beta(\omega_0,\Aoo,\Aot) < \sqrt{2}$, {\it i.e.} one has 
$\beta(\Aoo,\Aot) < \sqrt{2}$. (Note that in {\it massless} theories, such 
strictly spacelike separated wedge algebras {\it are}, in fact, maximally 
correlated \cite{34}.)  

     We remark that a consequence of Corollary 3.3 is that for a 
norm dense set of vector states, the maximal Bell correlation across 
two strictly spacelike separated wedge algebras lies strictly 
between 1 and $\sqrt{2}$, at least for suitably small (but nonzero) 
spacelike separation \cite{39}. But it is still {\it a priori} possible
that $\beta(\As(W_a),\As(W')) = 1$. This possibility has been 
eliminated, as follows. First of all, the following lemma was proven,
which gives useful sufficient conditions under which the Bell correlation 
invariant coincides with the maximal Bell correlation in a particular
normal state.

\proclaim{Lemma 3.4} \cite{39} 
Let $(\ABs)$ be a pair of commuting von Neumann algebras on a Hilbert space 
$\Hs$. Consider an automorphism $\alpha \in \Aut(\Bh)$ such that 
$\alpha(\As) = \As$ and $\alpha(\Bs) = \Bs$ and such that for all 
$A \in \tilde{\As} \vee \tilde{\Bs}$ (where $\tilde{\As}$ is dense in $\As$ 
and $\tilde{\Bs}$ is dense in $\Bs$) one has 
$\alpha^n (A) \to \omega_0(A)\cdot \idty$ in the 
$\sigma$-weak topology of $\Bh$ as $n \to \infty$ for some normal state 
$\omega_0 \in \Bh_*$. Then 
$$\beta_*(\ABs) = \beta(\ABs) = \beta(\omega_0,\ABs)
\, .$$ 
\endproclaim
 
     This yields the following corollary in algebraic quantum field theory.
We shall call a spacelike cylinder any open spacetime region $\Os$ 
such that for some spacelike direction $\vec{a}$ it is true that 
$\Os = \Os + t\vec{a}$ for all $t \in \RR$. Of course, such regions are 
necessarily unbounded, and wedges are examples of spacelike cylinders. 
 
\proclaim{Corollary 3.5} \cite{39} Let $\Os_1$, $\Os_2$ be parallel 
spacelike cylinders in three or more spacetime dimensions with 
$\Os_1 \subset \Os_2'$. Then in any
irreducible vacuum representation $\vacrept$ of a local net $\net$ 
one has the equality 
$$\beta(\Aoo,\Aot) = \beta(\omega_0,\Aoo,\Aot) \, ,
$$ 
where $\omega_0$ is the (unique) normal vacuum state on the 
representation. 
\endproclaim
 
     Taken in conjunction with the vacuum cluster bound, Corollary 3.3, this 
entails the following theorem. 
 
\proclaim{Theorem 3.6} \cite{39} Consider the function 
$$ b(a)=\beta(\Awo,\Awt) 
\, , $$
where $a\geq0$ is the spacelike separation of the wedges $W_1,W_2$. 
In an irreducible vacuum representation with a mass gap $m>0$, the function 
$b$ is lower semicontinuous and nonincreasing, and  $1\leq b(a)\leq \sqrt2$. 
The upper bound is an equality if and only if $a=0$. Consequently, 
the quantity $\beta(\Awo,\Awt)$ takes infinitely many different 
values (depending on the spacelike separation of the wedges), with an
accumulation point at $\sqrt2$. 
\endproclaim

\nind Hence, the invariant $\beta\pair$ can distinguish between infinitely
many different isomorphism classes of pairs $\pair$. At this point it is not 
known whether $\beta(\cdot,\cdot)$ can take any value between $1$ and $\sqrt2$.
\par

\bigpagebreak

\proclaim{Acknowledgements} This paper grew out of a very pleasant 
collaboration with Reinhard Werner. Part of this manuscript was written while 
the author was visiting the University of Rome under a grant from the CNR in 
the summer of 1996. For this opportunity he would like to thank his many
colleagues at the University of Rome, particularly Profs. Claudio D'Antoni,
John E. Roberts, and Laszlo Zsido.  \endproclaim

\bigpagebreak

\heading
References
\endheading

\eightpoint

\item{[1]}H. Araki, {\it Asymptotic ratio set and property $L_{\lambda}'$}, 
Publ. RIMS, Kyoto Univ., {\bf 6} (1970-1971), 443-460.

\item{[2]}H. Araki and E.J. Woods, {\it A classification of factors}, 
Publ. RIMS, Kyoto Univ., {\bf 4} (1968), 51-130.

\item{[3]}J.S. Bell, {\it On the Einstein-Podolsky-Rosen paradox}, Physics, 
{\bf 1} (1964), 195-200.

\item{[4]}J.S. Bell, {\it On the problem of hidden variables in quantum 
mechanics}, Reviews of Modern Physics, {\bf 38} (1966), 447-452. 

\item{[5]}N. Bohr, {\it Can quantum-mechanical description of reality be 
considered complete?}, Phys. Rev., {\bf 48} (1935), 696-702.  

\item{[6]}H. Br\'ezis, L. Nirenberg and G. Stampacchia, {\it A Remark on Ky 
Fan's Minimax principle}, Bolletino U.M.I., {\bf 6} (1972), 293-300.

\item{[7]}D. Buchholz, {\it Product states for local algebras}, Commun. Math. 
Phys., {\bf 36} (1974), 287-304.

\item{[8]}D. Buchholz and K. Fredenhagen, {\it Locality and the structure of 
particle states}, Commun. Math. Phys., {\bf 84} (1982), 1-54.

\item{[9]}D. Buchholz and E.H. Wichmann, {\it Causal independence and the 
energy-level density of states in local quantum field theory}, Commun. Math. 
Phys., {\bf 106} (1986), 321-344.

\item{[10]}D. Buchholz, C. D'Antoni and K. Fredenhagen, {\it The universal 
structure of local algebras}, Commun. Math. Phys., {\bf 111} (1987), 123-135.

\item{[11]}D. Buchholz and S.J. Summers, {\it An algebraic characterization of
vacuum states in Minkowski space}, Commun. Math. Phys., {\bf 155} (1993),
449-458.

\item{[12]}B.S. Cirel'son (Tsirelson), {\it Quantum generalizations of Bell's 
inequalities}, Lett. Math. Phys., {\bf 4} (1980), 93-100. 

\item{[13]}J.F. Clauser and M.A. Horne, {\it Experimental consequences of 
objective local theories}, Phys. Rev. {\bf D 10} (1974), 526-535. See also:
J.F. Clauser, M.A. Horne, A. Shimony and R.A. Holt, {\it Proposed experiment
to test local hidden-variable theories}, Phys. Rev. Lett., {\bf 23} (1969),
880-884.

\item{[14]}J.F. Clauser and A. Shimony, {\it Bell's theorem: Experimental tests
and implications}, Rep. Prog. Phys., {\bf 41} (1978), 1881-1927.

\item{[15]}W. DeBaere, {\it Einstein-Podolsky-Rosen paradox and Bell's 
inequalities}, Advances in Electronics, {\bf 68} (1986), 245-336.   

\item{[16]}S. Doplicher, R. Haag and J.E. Roberts, {\it Fields, observables and
gauge transformations, I, II}, Commun. Math. Phys., {\bf 13} (1969), 1-23, and
{\bf 15} (1969), 173-200.

\item{[17]}W. Driessler, {\it Duality and absence of locally generated 
superselection sectors for CCR-type algebras}, Commun. Math. 
Phys., {\bf 70} (1979), 213-220.

\item{[18]}A. Einstein, B. Podolsky and N. Rosen, Can quantum-mechanical 
description of physical reality be considered complete?, Physical Review, 
{\bf 47} (1935), 777-780.   

\item{[19]}M. Florig and S.J. Summers, {\it On the statistical independence of 
algebras of observables}, to appear in J. Math. Phys. 

\item{[20]}K. Fredenhagen, {\it A remark on the cluster theorem}, Commun. Math.
Phys., {\bf 97} (1985), 461-463.

\item{[21]}R. Haag, {\it Local Quantum Physics}, Springer Verlag, Berlin 
(1992).

\item{[22]}L.J. Landau, {\it Experimental tests of general quantum theories},
Lett. Math. Phys., {\bf 14} (1987), 33-40.

\item{[23]}L.J. Landau, {\it Experimental tests for distributivity}, Lett. 
Math. Phys., {\bf 25} (1992), 47-50.   

\item{[24]}G. Ludwig, {\it An Axiomatic Basis for Quantum Mechanics, I}, 
Springer, New York (1985).

\item{[25]}J. von Neumann, {\it Collected Works}, Volume {\bf 3}, Pergamon 
Press, New York, Oxford and London (1961).

\item{[26]}I. Pitowsky, {\it Quantum Probability - Quantum Logic}, 
Springer-Verlag, Berlin and New York (1989).  

\item{[27]}S. Popa, {\it Classification of Subfactors and Their Endomorphisms},
American Mathematical Society, Providence (1995).

\item{[28]}R.F. Powers, {\it Representations of uniformly hyperfinite 
algebras and their associated von Neumann rings}, Ann. Math., 
{\bf 86} (1967), 138-171.

\item{[29]}R.F. Powers, {\it UHF algebras and their applications to 
representations of the anticommutation relations}, in: {\it Carg\`ese 
Lectures in Physics}, Gordon and Breach, New York (1970).

\item{[30]}V. Shulman, private communication.

\item{[31]}S.J. Summers, {\it Normal product states for Fermions and twisted 
duality for CCR- and CAR-type algebras with application to the 
Yukawa$_2$ quantum field model}, Commun. Math. Phys., {\bf 86} (1982), 
111-141.

\item{[32]}S.J. Summers and R.F. Werner, {\it The vacuum violates Bell's 
inequalities}, Phys. Lett., {\bf A 110} (1985), 257-259.  

\item{[33]}S.J. Summers and R.F. Werner, {\it Bell's inequalities and quantum 
field theory, I: General setting}, preprint, unabridged version of 
\cite{34}, available from authors.  

\item{[34]}S.J. Summers and R.F. Werner, {\it Bell's inequalities and quantum 
field theory, I: General setting}, J. Math. Phys., {\bf 28} (1987), 2440-2447.
     
\item{[35]}S.J. Summers and R.F. Werner, {\it Bell's inequalities and quantum 
field theory, II: Bell's inequalities are maximally violated in the 
vacuum}, J. Math. Phys., {\bf 28} (1987), 2448-2456.     

\item{[36]}S.J. Summers and R.F. Werner, {\it Maximal violation of Bell's 
inequalities is generic in quantum field theory}, Commun. Math. 
Phys., {\bf 110} (1987), 247-259.  

\item{[37]}S.J. Summers and R.F. Werner, {\it Maximal violation of Bell's 
inequalities for algebras of observables in tangent spacetime 
regions}, Ann. Inst. Henri Poincar\'e, {\bf 49} (1988), 215-243.  

\item{[38]}S.J. Summers, {\it On the independence of local algebras in quantum 
field theory}, Rev. Math. Phys., {\bf 2} (1990), 201-247.   

\item{[39]}S.J. Summers and R.F. Werner, {\it On Bell's inequalities and 
algebraic invariants}, Lett. Math. Phys., {\bf 33} (1995), 321-334.

\item{[40]}A.M. Vershik and B.S. Tsirelson, {\it Formulation of Bell type 
problems and ``noncommutative'' convex geometry}, Adv. Sov. Math., {\bf 9} 
(1992), 95-114. 

\item{[41]}R.F. Werner, {\it Bell's inequalities and the reduction of 
statistical theories}, in: {\it Reduction in Science}, edited by W. Balzer, 
{\it et alia}, D. Reidel, Amsterdam (1984).

\item{[42]}R.F. Werner, {\it Remarks on a quantum state extension problem}, 
Lett. Math. Phys., {\bf 19} (1990), 319-326.

\enddocument
\bye